\documentclass[10pt,a4paper]{article}
\usepackage{amsmath,amssymb,amsthm}
\usepackage[T1]{fontenc}
\usepackage[latin1]{inputenc}

\theoremstyle{plain}
\newtheorem{theorem}{Theorem}[section]
\newtheorem{remark}{Remark}[section]
\newtheorem{corollary}{Corollary}[section]
\newtheorem{definition}{Definition}[section]
\newtheorem{proposition}{Proposition}[section]
\newtheorem{lemma}{Lemma}[section]

\usepackage{hyperref}

\newcommand{\F}{\mathcal{F}}
\renewcommand{\P}{\mathbb{P}}
\newcommand{\R}{\mathbb{R}}

\newcommand{\I}{1{\hskip -2.5 pt}\hbox{I}}
\newcommand{\Q}{\mathcal{Q}}
\newcommand{\C}{\mathcal{C}}
\newcommand{\D}{\mathcal{D}}

\newtheorem{property}{Property}
\def \follmer {F\"{o}llmer}
\def \ito {It\^{o} }
\def \cadlag {c\`adl\`ag }

\begin{document}

\title{On pathwise quadratic variation \\for c\`adl\`ag functions}
\author{Henry CHIU   \footnote{Department of Mathematics, Imperial College London. Email: {\tt h.chiu16@imperial.ac.uk}
Supported by  EPSRC Doctoral Training grant  1824430.} \ \& \ Rama CONT\footnote{Laboratoire de Probabilit\'es, Statistiques et Mod\'elisation, CNRS-Sorbonne Universit\'e.
Email: {\tt Rama.Cont@math.cnrs.fr}}
}
\date{June 2018. Final revision: October 2018.\\
Published in: \\Electronic Communications in Probability.}

\maketitle

\begin{abstract}
We revisit \follmer's concept of  quadratic variation of a \cadlag function along a sequence of time partitions  and discuss its relation with the Skorokhod topology. We show that in order to obtain a robust notion of pathwise quadratic variation applicable to sample paths of \cadlag processes, one must reformulate the definition of pathwise quadratic variation as a limit in Skorokhod topology of discrete approximations along the partition. 
One then obtains a simpler definition which implies the Lebesgue decomposition of the pathwise quadratic variation as a result, rather than requiring it as an extra condition. 
\end{abstract}
\tableofcontents

\section{Quadratic variation along a sequence of partitions}

 In his seminal paper {\it Calcul d'\ito sans probabilit\'es} \cite{follmer1981}, Hans \follmer{}  introduced  a pathwise concept of quadratic variation and used it to provide a pathwise proof of the \ito formula.
 F\"ollmer showed that if a function 
 $x\in D([0,T],\mathbb{R})$ has { quadratic variation} along a sequence  $\pi^n=(t^n_0=0< ..<t^n_j<...<t^n_{m(n)}=T)$ of time partitions of $[0,T]$ in the sense that for each $t\in [0,T]$ the limit
  \begin{equation}
		\mu^n := \sum_{[t^n_j, t^n_{j + 1}] \in \pi_n} \delta (\cdot - t_j)\quad |x(t^n_{j + 1}) - x(t^n_j)|^2 \label{eq.mun}
	\end{equation}
	converges weakly to a Radon measure $\mu_\pi$ such that\begin{center} $[x]^c_\pi: t\mapsto \mu_\pi([0,t]) -\sum_{s\leq t}|\Delta x(s)|^2$ is a continuous increasing function, \ \ \ (L)\end{center}
then a pathwise \ito formula may be obtained for functions of $x$ \cite{follmer1981}: for any $f\in C^2(\mathbb{R},\mathbb{R})$, 
 \begin{eqnarray}
		f(x(t))-f(x(0)) =& \int_0^t f'(x)d^\pi x + \frac{1}{2} \int_0^t f''(x(s))d[x]^c_\pi(s) \nonumber\\
		+&\mathop{\sum}_{0\leq s\leq t} \quad\left(\ f(x(s))-f(x(s-))-f'(x(s-))\Delta x(s) \ \right)\label{eq.ito}
	\end{eqnarray}
    where $\int_0^t f'(x)d^\pi x$ is a pathwise integral defined as a limit of left Riemann sums computed along the partition:
    \begin{equation}\label{eq.pathwiseintegral}
    \int_0^t g(x)d^\pi x:=\mathop{\lim}_{n\to\infty} \sum_{\pi_n} g(x(t^n_j)).\left( x(t^n_{j + 1}) - x(t^n_j)\right).
    \end{equation}
The quantity 
$$[x]_\pi(t)= \mu_\pi([0,t])=[x]^c_\pi(t)+ \sum_{s\leq t}|\Delta x(s)|^2$$
is called the quadratic variation of $x$ along  $\pi$.
 This result has many interesting applications  and has been  extended  to less regular functions \cite{bertoin1987,cp2018,davis2018,perkowski2015} and path-dependent functionals  \cite{ananova2017,cont2012,CF09,CF10B,cp2018}. With the exception of \cite{CF10B,vovk2015,perkowski2017}, these extensions have focused on continuous paths.

\follmer's definition \cite{follmer1981} contains the condition (L) on the Lebesgue decomposition of the limit $\mu_\pi$: the atoms of $\mu$ should correspond exactly to the jumps of $x$ and their mass should be $|\Delta x(t)|^2$ or, equivalently, the discontinuity points of $[x]_\pi$ should coincide with those of $x$, with $\Delta [x]_\pi(t)=|\Delta x(t)|^2$.
This condition can not be removed: as shown by Coquet et al. \cite{coquet2006}, there are counterexamples of continuous functions $x$ such that \eqref{eq.mun} converges to a limit with atoms. Conversely, one can give examples of discontinuous functions for which \eqref{eq.mun} converges to an atomless measure.
	If this condition is not satisfied, then the pathwise  integral \eqref{eq.pathwiseintegral} fails to satisfy  \begin{eqnarray*}
    \Delta\left(\int_0^tx(s-)dx(s)\right)(t)= x(t-)\Delta x(t),
\end{eqnarray*} at each $t$ where condition (L) is not met.
On the other hand, this condition (L) is not easy to check and seems to require a link between the path $x$ and the sequence of partitions $\pi$, making it difficult to apply to sample paths of stochastic processes.
 
 In this work, we revisit \follmer's concept of pathwise quadratic variation along a sequence of partitions and show that it has hitherto unsuspected links with the Skorokhod topology. In particular, we show that in order to obtain a robust notion of pathwise quadratic variation applicable to sample paths of \cadlag processes, one must reformulate the definition of the quadratic variation  as a limit, in Skorokhod topology, of discrete approximations defined along the partition. 
This leads to a simpler definition of pathwise quadratic variation which holds in any dimension and, rather than requiring the Lebesgue decomposition of the pathwise quadratic variation as an extra condition, yields it as a consequence.

{\bf Outline} We begin by recalling  \follmer's definition of pathwise quadratic variation and variations of it which have been used in the literature. We then introduce a new definition of quadratic variation for real-valued \cadlag functions based on the Skorokhod topology and prove equivalence among the various definitions. Section \ref{sec.multidimensional} extends the results to vector-valued functions: we show that, unlike \follmer's original definition in which the one dimensional case plays a special role, our definition applies regardless of dimension, thus simplifying various statements regarding quadratic variation for vector-valued functions.
Finally, in Section \ref{sec.semimartingale}, we show that our approach leads to simple proofs for various properties of pathwise quadratic variation.

\section{Pathwise quadratic variation for cadlag functions}
\label{sec.skorokhod}
Let $\pi:=(\pi_n)_{n\geq 1}$ be a sequence of partitions $\pi_n=(t^{n}_0,...,t^{n}_{k_n})$ of $[0,\infty)$ into intervals $0=t^{n}_0<...<t^{n}_{k_n}<\infty$; $t^{n}_{k_n}\uparrow\infty$ with vanishing mesh $|\pi_n|\downarrow 0$ on compacts. By convention, $\max(\emptyset\cap\pi_n):=0$, $\min(\emptyset\cap\pi_n):=t^{n}_{k_n}$.

Denote $\D:=\D([0,\infty),\R)$ the  space of cadlag functions and $\C:=\C([0,\infty),\R)$ the subspace of real-valued continuous functions. We equip $\D$ with a metric $d$ which induces the Skorokhod  $J_1$ topology 
\cite[Ch. VI]{jacodshiryaev}. Denote $\D^{+}_{0}\subset\D$ to be the subset of non-negative increasing right-continuous functions null at $0$.

\begin{definition}[\follmer 1981]\label{def:2.1}
 $x\in \D$ has (finite) quadratic variation $[x]_\pi$ along $\pi$ if the sequence of  measures \begin{eqnarray}\label{eq:2.1}\mu_n:=\sum_{ \pi_n}(x(t^n_{i+1})-x(t^n_{i}))^2 \delta(t^n_{i})\end{eqnarray} converges vaguely to a Radon measure $\mu$ on $[0,\infty)$ with $[x]_\pi(t)=\mu([0,t])$, such that $[x]^c_\pi$ defined by $[x]_\pi^c(t):=\mu([0,t])- \sum_{s\leq t}(\Delta x_s)^2$ is a continuous increasing function.
$[x]_\pi(t)=\mu([0,t])$ is then called the quadratic variation of $x$ along $\pi$, and admits the following Lebesgue decomposition:
\begin{eqnarray}\label{eq:2.2}
[x]_\pi(t)=[x]_\pi^{c}(t)+\sum_{s\leq t}|\Delta x(s)|^2.
\end{eqnarray}
We denote $\Q_0^{\pi}$ the set of $x\in \D$ satisfying these properties.
\end{definition}

	At this point let us  point out a link between vague and weak convergence of Radon measures on $[0,\infty)$, a link which is well known in the case where the measures are sub-probability measures:
\begin{lemma}\label{lem:2.2}
Let $v_n$ and $v$ be non-negative Radon measures on $[0,\infty)$ and $J\subset[0,\infty)$ be the set of atoms of $v$, the followings are equivalent:\\

\noindent(i) $v_n\rightarrow v$ vaguely on $[0,\infty)$.\\
(ii) $v_n\rightarrow v$ weakly on $[0,T]$ for every $T\notin J$.
\end{lemma}
\begin{proof}
Let $f\in C_K([0,\infty))$ be a compactly supported continuous function. Since $J$ is countable, $\exists$ $T\notin J$; supp$(f)\subset[0,T]$. Now (ii) $\Rightarrow\int_{0}^{\infty} f dv_n=\int_{0}^{T}f dv_n \longrightarrow \int_{0}^{T}f dv=\int_{0}^{\infty} f dv\Rightarrow$ (i). Suppose (i) holds, let $T\notin J$ and $f\in \C\left([0,T],\|\cdot\|_{\infty}\right)$. Since $f=(f)^{+}-(f)^{-}$, we may take $f\geq 0$ and define the following extensions:\begin{eqnarray*}\overline{f}^{\epsilon}(t)&:=&f(t)\I_{[0,T]}(t)+f(T)\left(1+\frac{T-t}{\epsilon}\right)\I_{(T,T+\epsilon]}(t)\\
\underline{f}^{\epsilon}(t)&:=&f(t)\I_{[0,T-\epsilon]}(t)+f(T)\left(\frac{T-t}{\epsilon}\right)\I_{(T-\epsilon,T]}(t)
,\end{eqnarray*}then $\overline{f}^{\epsilon}$, $\underline{f}^{\epsilon}\in\C_K([0,\infty))$, $0\leq\underline{f}^{\epsilon}\leq f\I_{[0,T]}\leq\overline{f}^{\epsilon}\leq\|f\|_{\infty}$ and we have \begin{eqnarray*}\int_{0}^{\infty}\underline{f}^{\epsilon}dv_{n}\leq\int_{0}^{T} f dv_{n}\leq\int_{0}^{\infty}\overline{f}^{\epsilon}dv_n.\end{eqnarray*}Since $v_n \rightarrow v$ vaguely and $T\notin J$, thus, as $n\rightarrow\infty$ we obtain \begin{eqnarray*}0&\leq&\limsup_{n}\int_{0}^{T} f dv_{n}-\liminf_{n}\int_{0}^{T} f dv_{n}\leq\int_{0}^{\infty}\overline{f}^{\epsilon}-\underline{f}^{\epsilon}dv\\
&\leq&\|f\|_{\infty}v\left((T-\epsilon,T+\epsilon]\right)\stackrel{\epsilon}{\longrightarrow}0,\end{eqnarray*}hence by monotone convergence\begin{eqnarray*}\lim_{n}\int_{0}^{T} f dv_{n}=\lim_{\epsilon}\int_{0}^{\infty}\underline{f}^{\epsilon}dv=\int_{0}^{T}fdv\end{eqnarray*}and (ii) follows.
\end{proof}\noindent

\begin{proposition}\label{prop:2.3}
If $x\in \Q_0^{\pi}$, then the pointwise limit $s$ of \begin{eqnarray}\label{eq:2.3}s_n(t):=\sum_{t_i\in\pi_{n}}|x(t_{i+1}\wedge t)-x(t_{i}\wedge t)|^2\end{eqnarray} exists, $s=[x]_\pi$ and $s$ admits the Lebesgue decomposition:\begin{eqnarray}\label{eq:2.4}
s(t)=s^{c}(t)+\sum_{s\leq t}(\Delta x_s)^2.
\end{eqnarray}
\end{proposition}
\begin{proof}
If $x\in Q_{0}^{\pi}$, define \begin{eqnarray*}q_n(t):=\sum_{\pi_n\ni t_i\leq t}(x(t_{i+1})-x(t_{i}))^2,\end{eqnarray*}the distribution function of $\mu_n$ in (\ref{eq:2.1}). Since $\mu_n\rightarrow \mu$ vaguely, we have $q_n\rightarrow [x]$ pointwise at all continuity points of $[x]$ (Lemma \ref{lem:2.2} \& \cite[X.11]{doob}). Let $I$ be the set of continuity points of $[x]$. Observe $(q_n)$ is monotonic in $[0,\infty)$ and $I$ is dense in $[0,\infty)$, if $t\notin I$, it follows \cite[X.8]{doob} that \begin{eqnarray*}[x]({t-})\leq\liminf_{n} q_{n}(t)\leq\limsup_{n} q_{n}(t)\leq[x]({t+})=[x](t).\end{eqnarray*}Thus, we may take any subsequence $(n_k)$ such that $\lim_{k} q_{n_k}(t)=:q(t)$. Since $x\in Q_0^{\pi}$ and the Lebesgue decomposition (\ref{eq:2.2}) holds on $[x]$, we have \begin{eqnarray}\label{eq:2.5}\stackrel{\geq 0}{\left([x](t+\epsilon)-q(t)\right)}+\stackrel{\geq 0}{\left(q(t)-[x](t-\epsilon)\right)}=[x](t+\epsilon)-[x]({t-\epsilon})\mathop{\to}^{\epsilon\to 0}|\Delta x(t)|^2.\end{eqnarray}If $t\pm\epsilon\in I$, $\tilde{\pi}_k:=\pi_{n_k}$ and $t_j^{(k)}:=\max\{\tilde{\pi}_k\cap[0,t)\}$, the second sum in (\ref{eq:2.5}) is 
\begin{eqnarray*}
\lim_{k}\sum_{\substack{t_i\in\tilde{\pi}_k; \\t-\epsilon< t_i\leq t}}(x(t_{i+1})-x(t_{i}))^2
=\lim_{k}\sum_{\substack{t_i\in\tilde{\pi}_k; \\t-\epsilon< t_i<t_j^{(k)} }}\stackrel{\geq 0}{(x(t_{i+1})-x(t_{i})^2}+(\Delta x(t))^2\geq(\Delta x(t))^2
\end{eqnarray*}by the fact that $x$ is \cadlag and that $t\notin I$. 

	We see  from (\ref{eq:2.5}) that $q(t)=[x](t)$ as $\epsilon\rightarrow 0$. Since the choice of the convergent subsequence is arbitrary, we conclude that $q_n\rightarrow[x]$ pointwise on $[0,\infty)$. Observe that the pointwise limits of $(s_n)$ and $(q_n)$ coincide i.e. \begin{eqnarray}|s_n(t)-q_n(t)|=(x_{t^{(n)}_{i+1}}-x(t))^2+2(x_{t^{(n)}_{i+1}}-x(t))(x(t)-x_{t^{(n)}_{i}})\label{eq:2.6}\end{eqnarray}converges to $0$ by the right-continuity of $x$, where $t_{i}^{(n)}:=\max\left\{\pi_n\cap[0,t]\right\}$ and that $x\in Q_{0}^{\pi}$, Prop.~\ref{prop:2.3} follows.
\end{proof}\noindent\\

	Denote $\Q_1^{\pi}$  the set of  $x\in\D$ such that   $(s_n)$ defined in \eqref{eq:2.3} has a pointwise limit $s$ with  Lebesgue decomposition given by \eqref{eq:2.4}.
	Then $\Q_0^{\pi}\subset\Q_1^{\pi}$ and we have:
\begin{proposition}\label{prop:2.4}
If $x\in \Q_1^{\pi}$, then the pointwise limit $q$ of \begin{eqnarray}\label{eq:2.7}q_n(t):=\sum_{\pi_n\ni t_i\leq t}(x(t_{i+1})-x(t_{i}))^2\end{eqnarray} exists. Furthermore $q=s$ and admits the Lebesgue decomposition:\begin{eqnarray}\label{eq:2.8}
q(t)=q^{c}(t)+\sum_{s\leq t}(\Delta x_s)^2.
\end{eqnarray}
\end{proposition}
\begin{proof}
Since the pointwise limits of $(s_n)$ and $(q_n)$ coincide by (\ref{eq:2.6}). Prop.~\ref{prop:2.4} now follows from $x\in \Q_1^{\pi}$.
\end{proof}

\noindent
	Denote $\Q_2^{\pi}$  the set of  $x\in\D$ such that the quadratic sums $(q_n)$ defined in \eqref{eq:2.7}
	have a  pointwise limit $q$ with Lebesgue decomposition  \eqref{eq:2.8}.
	Then $\Q_0^{\pi}\subset\Q_1^{\pi}\subset\Q_2^{\pi}$ and we have:
\begin{proposition}\label{prop:2.5}
If $x\in\Q_2^{\pi}$, then $q_n\rightarrow q$ in the Skorokhod topology. 
\end{proposition}
\begin{proof}
Since $x\in\Q_2^{\pi}$, we have $q_n\rightarrow q$ pointwise on $[0,\infty)$ and that $(q_n)$, $q$ are elements in $\D^{+}_{0}$. By \cite[Thm.VI.2.15]{jacodshiryaev}, it remains to show that \begin{eqnarray*}\sum_{s\leq t}(\Delta q_n(s))^2\stackrel{n\to\infty}{\longrightarrow}\sum_{s\leq t}(\Delta q(s))^2\end{eqnarray*} on a dense subset of $[0,\infty)$. Let $t>0$, define $J^{\epsilon}:=\{s\geq 0| (\Delta X_s)^2\geq\frac{\epsilon}{2}\}$, $J^{\epsilon}_{n}:=\{t_i\in\pi_n | \exists s\in(t_i,t_{i+1}]; (\Delta X_s)^2\geq\frac{\epsilon}{2}\}\subset\pi_n$ and observe that \begin{eqnarray}\label{eq:2.9}
\sum_{s\leq t}(\Delta q_n(s))^2&=&\sum_{\pi_n\ni t_i\leq t}(x(t_{i+1})-x(t_{i}))^4\nonumber\\
&=&\sum_{J^{\epsilon}_{n}\ni t_i\leq t}(x(t_{i+1})-x(t_{i}))^4+\sum_{(J^{\epsilon}_{n})^{c}\ni t_i\leq t}(x(t_{i+1})-x(t_{i}))^4.
\end{eqnarray}Since $x$ is \cadlag and that $|\pi_n|\downarrow 0$ on compacts, the first sum in (\ref{eq:2.9}) converges to $\sum_{J^{\epsilon}\ni s\leq t}(\Delta x_s)^4$ and the second sum in (\ref{eq:2.9}) \begin{eqnarray*}\sum_{(J^{\epsilon}_{n})^{c}\ni t_i\leq t}(x(t_{i+1})-x(t_{i}))^4\leq\left(\sup_{(J^{\epsilon}_{n})^{c}\ni t_i\leq t}(x(t_{i+1})-x(t_{i}))^2\right)\ \sum_{(J^{\epsilon}_{n})^{c}\ni t_i\leq t}(x(t_{i+1})-x(t_{i}))^2\leq\epsilon q(t)\end{eqnarray*}for sufficiently large $n$ \cite[Appendix A.8]{CF10B} hence \begin{eqnarray*}
\lim_{n}\sum_{s\leq t}(\Delta q_n(s))^2=\sum_{J^{\epsilon}\ni s\leq t}(\Delta x_s)^4+\overbrace{\limsup_{n}{\sum_{(J^{\epsilon}_{n})^{c}\ni t_i\leq t}(x(t_{i+1})-x(t_{i}))^4} }^{\leq\epsilon q(t)}.
\end{eqnarray*}
By the Lebesgue decomposition (\ref{eq:2.8}), we observe $\sum_{J^{\epsilon}\ni s\leq t}(\Delta x_s)^4\leq q(t)^2$ and that\begin{eqnarray*}
\lim_{n}\sum_{s\leq t}(\Delta q_n(s))^2=\sum_{s\leq t}(\Delta x_s)^4=\sum_{s\leq t}(\Delta q(s))^2\end{eqnarray*}as $\epsilon\rightarrow 0$.                         
\end{proof}\noindent\\
	Denote $\Q^{\pi}$ the set of \cadlag functions $x\in\D$ such that the limit $\tilde{q}$ of $(q_n)$ exists in $(\D,d)$. Then $\Q_0^{\pi}\subset\Q_1^{\pi}\subset\Q_2^{\pi}\subset\Q^{\pi}$ and we have:
\begin{proposition}\label{prop:2.6}
$\Q^{\pi}\subset\Q_{0}^{\pi}$ and $\tilde{q}=[x]$.
\end{proposition}
\begin{proof}
Let $x\in\Q^{\pi}$ and $I$ be the set of continuity points of $\tilde{q}$. \cite[VI.2.1(b.5)]{jacodshiryaev} implies that $q_n\rightarrow\tilde{q}$ pointwise on $I$. Since $q_n\in\D^{+}_{0}$ and $I$ is dense on $[0,\infty)$, it follows $\tilde{q}\in\D^{+}_{0}$. Denote $\mu$ to be the Radon measure of $\tilde{q}$ on $[0,\infty)$, observe the set of atoms of $\mu$ is $J:=[0,\infty)\backslash I$ and that $(q_n)$ are the distribution functions of the discrete measures $(\mu_n)$ in (\ref{eq:2.1}). Thus, by (Lemma \ref{lem:2.2} \& \cite[X.11]{doob}), we see that $\mu_n\longrightarrow\mu$ vaguely on $[0,\infty)$. 
	
	If $t>0$, put $t^{(n)}_i:=\max\{\pi_n\cap[0,t)\}$. Since $|\pi_n|\downarrow 0$ on compacts, we have $t^{(n)}_i<t$, $t^{(n)}_i\uparrow t$ and $t^{(n)}_{i+1}\downarrow t$. Observe that \begin{eqnarray}\label{eq:2.10}\Delta q_n(t)=\begin{cases}
    \left(x(t_{i+1})-x(t_i)\right)^2, & \text{if $t=t_i\in\pi_n$}.\\
    0, & \text{otherwise}.
  \end{cases}
\end{eqnarray}If $\Delta\tilde{q}(t)=0$, \cite[VI.2.1(b.5)]{jacodshiryaev} implies that $\Delta q_n(t^{(n)}_i)\rightarrow \Delta\tilde{q}(t)$. Hence, by the fact that $x$ is \cadlag, $(\Delta x(t))^2=\lim_{n}\Delta q_n(t^{(n)}_i)=\Delta\tilde{q}(t)$. If $\Delta\tilde{q}(t)>0$, there exists \cite[VI.2.1(a)]{jacodshiryaev} a sequence $t'_n\rightarrow t$ such that $\Delta q_n(t'_n)\rightarrow\Delta\tilde{q}(t)>0$. Using the fact that $x$ is \cadlag, $t'_n\rightarrow t$, (\ref{eq:2.10}) and \cite[VI.2.1(b)]{jacodshiryaev}, we deduce that $(t'_n)$ must coincide with $(t^{(n)}_i)$ for all $n$ sufficiently large, else we will contradict $\Delta\tilde{q}(t)>0$. Thus, $(\Delta x(t))^2=\lim_{n}\Delta q_n(t^{(n)}_i)=\lim_{n}\Delta q_n(t'_n)=\Delta\tilde{q}(t)$ and the Lebesgue decomposition (\ref{eq:2.2}) holds on $\tilde{q}$.
 
	By Def.\ref{def:2.1}, we have $\tilde{q}=[x]$ hence $\Q^{\pi}\subset\Q_{0}^{\pi}$.
\end{proof}

\begin{theorem}\label{thm:2.7}
Let
\begin{itemize}
\item  $\Q_0^{\pi}$ be the set of  $x\in\D$ satisfying Definition \ref{def:2.1}.
\item  $\Q_1^{\pi}$  the set of  $x\in\D$ such that   $(s_n)$ defined in \eqref{eq:2.3} has a pointwise limit $s$ with  Lebesgue decomposition given by \eqref{eq:2.4}.
    \item $\Q_2^{\pi}$  the set of  $x\in\D$ such that the quadratic sums $(q_n)$ defined by \eqref{eq:2.6}
	have a  pointwise limit $q$ with Lebesgue decomposition  \eqref{eq:2.7}.
	\item $\Q^{\pi}$ the set of \cadlag functions $x\in\D$ such that the limit $\tilde{q}$ of $(q_n)$ exists in $(\D,d)$.
\end{itemize}
Then:
\begin{enumerate}
\item[(i)] $\Q_{0}^{\pi}=\Q_{1}^{\pi}=\Q_{2}^{\pi}=\Q^{\pi}$.
\item[(ii)] If $x\in\Q_{0}^{\pi}$, then $[x]_\pi=s(x)=q(x)=\tilde{q}(x)$.
\item[(iii)] $x$ has finite quadratic variation along $\pi$ if and only if \begin{eqnarray*}
q_n(t):=\sum_{\pi_n\ni t_i\leq t}(x(t_{i+1})-x(t_{i}))^2
\end{eqnarray*}converges in $(\D,d)$.
\item[(iv)] If $(q_n)$ defined by \eqref{eq:2.6} converges in $(\D,d)$, the limit is equal to $[x]$.
\end{enumerate}
\end{theorem}
    
\begin{proof}
These results are  a consequence of Prop.~\ref{prop:2.3}, \ref{prop:2.4}, \ref{prop:2.5} and \ref{prop:2.6}.
\end{proof}

	We see that the two defining properties of $[x]$ in Def.~\ref{def:2.1} are 
	consequences per Thm.~\ref{thm:2.7}. 
	The following corollary treats the special case of continuous functions.
\begin{corollary} 
Let $x\in\Q^{\pi}$, $s_n$  defined as in \eqref{eq:2.3}, $(q_n)$ defined by \eqref{eq:2.7}.\begin{enumerate}
\item[i] $q_{n}\rightarrow[x]$  uniformly on compacts in $[0,\infty)$ if and only if $x$ is continuous.
\item[ii] If $q_{n}\rightarrow[x]$  uniformly on compacts in $[0,\infty)$, then  
$s_{n}\rightarrow[x]$ uniformly on compacts.
\end{enumerate}
\end{corollary}	
\begin{proof}
(i): The if part follows from Prop.~\ref{prop:2.6}, (\ref{eq:2.2}) and \cite[VI.1.17(b)]{jacodshiryaev}. The only if part: 
Let $t\geq 0$, it is well known that \begin{eqnarray*}
\Delta q_{n}(t)\rightarrow \Delta [x](t)
\end{eqnarray*}by uniform convergence. Put $t'_n:=\max \{t_i<t|t_i\in \pi_n\}$, since $q_{n}\rightarrow[x]$ in the Skorokhod topology, we also have \begin{eqnarray*}
\Delta q_{n}(t'_n)\rightarrow \Delta [x](t)
\end{eqnarray*}by \eqref{eq:2.10} and \cite[VI.2.1(b)]{jacodshiryaev}. If $\Delta [x](t)>0$, then \cite[VI.2.1(b)]{jacodshiryaev} implies $t'_n$ must coincide with $t$ for all $n$ large enough, but $t'_n<t$ for all $n$, hence $\Delta [x](t)=0$ which implies $\Delta x(t)=0$ by Prop.~\ref{prop:2.6} and (\ref{eq:2.2}). Since $t$ is arbitrary, we conclude that $x\in\C$.\\

(ii): Let $T>0$, $\|\cdot\|(t)$ the supremum norm on $\D([0,T])$ and observe that \begin{eqnarray*}
\|s_{n}-[x]\|(t)\leq\|q_{n}-[x]\|(t)+\|s_{n}-q_{n}\|(t).
\end{eqnarray*}Since (i) implies $x\in\C$, (ii) now follows from uniform continuity of $x$ and (\ref{eq:2.8}). 
\end{proof}

\begin{remark}
 The converse of (ii) is  not true in general.
\end{remark}
\begin{remark}
We note that some references have used the pointwise limit of the sequence\begin{eqnarray*}
p_n(t):=\sum_{\pi_n\ni t_{i+1}\leq t}(x(t_{i+1})-x(t_{i}))^2
\end{eqnarray*} together with the Lebesgue decomposition (\ref{eq:2.2}), to define $[x]$. To see why this is not the correct choice, take  $t_0\notin\pi$, put $x(t):=\I_{[t_0,\infty)}(t)$ then obviously $[x](t_0)=\lim s_n(t_0)=\lim q_n(t_0)=1$ but $\lim p_n(t_0)=0$.
\end{remark}

\begin{remark}{\em
Vovk  \cite[Sec. 6]{vovk2015} proposes a different notion of pathwise quadratic variation along a sequence of partitions, which is shown to coincide with \follmer's definition under  the additional assumption that the sequence of partitions $\pi$ is refining  and exhausts all discontinuity points of the path  \cite[Prop. 6.3 \& 6.4]{vovk2015}. 

This requirement of exhausting all jumps can always be satisfied for a given \cadlag path by adding all discontinuity points to the sequence of partitions. However if one is interested in applying this definition to  a process, say a semimartingale, then in general there may exist no sequence of partitions satisfying this condition.
And, if this requirement of exhausting all discontinuity points is removed, then Vovk's definition differs from \follmer's (and therefore, fails to satisfy the Ito formula \eqref{eq.ito} in general). 

By contrast, our definition does not require such a condition and easily carries over to stochastic processes without requiring the use of random partitions (see Theorem \ref{thm:4.4}).}
\end{remark}

\section{Quadratic variation for multidimensional functions} 
\label{sec.multidimensional}

	Denote $\D^{m}:=\D([0,\infty),\R^{m})$ and $\D^{m\times m}:=\D([0,\infty),\R^{m\times m})$ to be the Skorokhod spaces \cite{billingsley1968,skorokhod1956,jakubowski1986}, each of which equipped with a metric $d$ which induces the corresponding Skorokhod   $J_1$ topology \cite[Ch. VI]{jacodshiryaev}. $\C^{m}:=\C([0,\infty),\R^{m})$ the subspace of continuous functions in $\D^{m}$. We recall Theorem \ref{thm:2.7} from the one dimensional case $n=1$ that (Def.~\ref{def:2.1}) and (Thm.~\ref{thm:2.7}.iii) are equivalent.
	
	It is known that  $x,y\in\Q^{\pi}$ does not imply $x+y \in\Q^{\pi}$ \cite{cont2016,schied2016} so one cannot for instance define a quadratic covariation $[x,y]_\pi$ of two such functions in the obvious way. This prevents a simple componentwise definition of the finite quadratic variation property for vector-valued functions. Therefore, the notion of quadratic variation in the multidimensional setting was originally defined in \cite{follmer1981} as follows:\\
	
\begin{definition}[\follmer 1981]\label{def:3.1}
We say that $\mathbf{x}:=(x^{1},\ldots,x^{m})^{T}\in \D^{m}$ has finite quadratic variation along $\pi$ if all $x^{i}$, $x^{i}+x^{j}$ $(1\leq i, j\leq m)$ have finite quadratic variation.\end{definition}
	The  quadratic (co)variation $[x^{i},x^{j}]$ is then defined as
\begin{eqnarray}\label{eq:3.1}
[x^{i},x^{j}]_\pi(t):=\frac{1}{2}\left([x^{i}+x^{j}]_\pi(t)-[x^{i}]_\pi(t)-[x^{j}]_\pi(t)\right),
\end{eqnarray} which admits the following Lebesgue decomposition:
\begin{eqnarray}\label{eq:3.2}
[x^{i},x^{j}]_\pi(t)=[x^{i},x^{j}]_\pi^{c}(t)+\sum_{s\leq t}\Delta x^{i}(s)\Delta x^{j}(s).
\end{eqnarray} The function $[\mathbf{x}]_\pi:=([x^{i},x^{j}])_{1\leq i, j\leq m}$, which takes values in the cone of symmetricsemidefinite  positive matrices, is called the quadratic (co)variation of $\mathbf{x}$.\\

	Note Def.~\ref{def:3.1} requires  first introducing the case $m=1$. The following definition, by contrast, avoids this and  directly defines the concept of  multidimensional quadratic variation in any dimension:
\begin{definition}\label{def:3.2}
 $\mathbf{x}\in \D^{m}$ has finite quadratic variation $[\mathbf{x}]_\pi$ along $\pi$ if
\begin{eqnarray*}
\mathbf{q}_n(t):=\sum_{\pi_n\ni t_i\leq t} (\mathbf{x}(t_{i+1})-\mathbf{x}(t_{i}) )(\mathbf{x}(t_{i+1})-\mathbf{x}(t_{i}) )^{T}
\end{eqnarray*} converges to $[\mathbf{x}]_\pi$ in $( \D^{m\times m},d)$. 
\end{definition}
We shall now prove the equivalence of these definitions.

	Define, for $u,v,w\in\D$ \begin{eqnarray*}
q_n^{(u,v)}(t):=\sum_{\pi_n\ni t_i\leq t}(u_{t_{i+1}}-u_{t_{i}})(v_{t_{i+1}}-v_{t_{i}})
\end{eqnarray*}and $q_n^{(w)}:=q_n^{(w,w)}$. Note that the Skorokhod topology on $(\D^m,d)$ is strictly finer than the product topology on $(\D,d)^{m}$ \cite[VI.1.21]{jacodshiryaev} and that $(\D,d)$ is not a topological vector space \cite[VI.1.22]{jacodshiryaev}. The following lemma is essential:
\begin{lemma}\label{lem:3.3}
Let $t>0$. There exists a sequence $t_{n}\rightarrow t$ such that for all $u,v\in\D$, if $(q_n^{(u,v)})$ converges in $(\D,d)$ then
\begin{eqnarray*}\lim_{n}\left(\Delta q_n^{(u,v)}(t_{n})\right)=\Delta \left(\lim_{n}q_n^{(u,v)}\right)(t).\end{eqnarray*} 
\end{lemma}
Note that the sequence $t_n$ is chosen from the partition points of $\pi_n$, independently of $u,v\in\D$.
\begin{proof}
Define $t^{(n)}_i:=\max\{\pi_n\cap[0,t)\}$. Since $|\pi_n|\downarrow 0$ on compacts, we have $t^{(n)}_i<t$, $t^{(n)}_i\uparrow t$ and $t^{(n)}_{i+1}\downarrow t$. Observe that \begin{eqnarray}\Delta q_n^{(u,v)}(t)=\begin{cases}\label{eq:3.3}
    \left(u_{t_{i+1}}-u_{t_i}\right)\left(v_{t_{i+1}}-v_{t_i}\right), & \text{if $t=t_i\in\pi_n$}.\\
    0, & \text{otherwise}.
  \end{cases}
\end{eqnarray}Put $\tilde{q}:=\lim_{n}q_n^{(u,v)}$. If $\Delta \tilde{q}(t)=0$, \cite[VI.2.1(b.5)]{jacodshiryaev} implies that $\Delta q_n^{(u,v)}(t^{(n)}_i)\rightarrow \Delta \tilde{q}(t)$. If $\Delta \tilde{q}(t)>0$, there exists \cite[VI.2.1(a)]{jacodshiryaev} a sequence $t'_n\rightarrow t$ such that $\Delta q^{(u,v)}_n(t'_n)\rightarrow\Delta \tilde{q}(t)>0$. Using the fact that $u,v$ are \cadlag, $t'_n\rightarrow t$ and (\ref{eq:3.3}), we deduce that $(t'_n)$ must coincide with $(t^{(n)}_i)$ for $n$ sufficiently large, else we will contradict $\Delta \tilde{q}(t)>0$. Put $t_n:=t^{(n)}_i$.
\end{proof}\noindent

\begin{proposition}\label{prop:3.4}
Let $x,y\in\Q^{\pi}$, then $(q_n^{(x+y)})$ converges in $(\D,d)$ if and only if $(q_n^{(x,y)})$ does. In this case, $x+y\in\Q^{\pi}$ and $\lim_n q_n^{(x,y)}=\frac{1}{2}\left([x+y]-[x]-[y]\right)$. 
\end{proposition}
\begin{proof}
Since \begin{eqnarray*}q_n^{(x+y)}=q_n^{(x)}+q_n^{(y)}+2 q_n^{(x,y)}\end{eqnarray*} and that $x,y\in\Q^{\pi}$, Prop.~\ref{prop:3.4} follows from Lemma~\ref{lem:3.3} and \cite[VI.2.2(a)]{jacodshiryaev}.  
\end{proof}

\begin{proposition}\label{prop:3.5}
$(\mathbf{q}_n)$ converges in $(\D^{m\times m},d)$ if and only if it converges in $(\D,d)^{m\times m}$.
\end{proposition}
\begin{proof}
Since the Skorokhod topology on $(\D^{m\times m},d)$ is strictly finer than the product topology on $(\D,d)^{n\times n}$ \cite[VI.1.21]{jacodshiryaev}, we have $(\D^{m\times m},d)$ convergence implies $(\D,d)^{m\times m}$ convergence. The other direction follows  from the observation that \begin{eqnarray}\label{eq:3.4}\mathbf{q}_n=\left(q_n^{(x^{i},x^{j})}\right)_{1\leq i\leq j\leq n},\end{eqnarray} satisfies Lemma~\ref{lem:3.3} and \cite[VI.2.2(b)]{jacodshiryaev}.
\end{proof}

\begin{theorem}\label{thm:3.6}
Definitions~\ref{def:3.1} and \ref{def:3.2} are equivalent.
\end{theorem}
\begin{proof}
This is a  consequence of (\ref{eq:3.1}), (\ref{eq:3.4}), Prop.~\ref{prop:3.4} \& \ref{prop:3.5} and Thm.~\ref{thm:2.7}
\end{proof}
	
\begin{corollary}\label{cor:3.7}
If $\mathbf{x}\in\D^{m}$ has finite quadratic variation, then\\

\noindent(i) $\mathbf{q_n}\rightarrow [\mathbf{x}]$ locally uniformly on $[0,\infty)$ if and only if $\mathbf{x}\in\C^{n}$.\\
(ii)  $F(\mathbf{q_n})\rightarrow F([\mathbf{x}])$ for all  functionals $F$ which are $J_1$-continuous at $[\mathbf{x}]$.\\
\end{corollary}
\begin{proof}
This is a consequence of Thm.\ref{thm:3.6}, (\ref{eq:3.2}) and \cite[VI.1.17.b]{jacodshiryaev}.
\end{proof}

\begin{remark}\label{rem:3.8}
For $\mathbf{x}$ to have finite quadratic variation, it is sufficient that $(\mathbf{q_n})$ converges in $(\D,d)^{m\times m}$ due to Prop.~\ref{prop:3.5}. (i.e. component-wise convergence) 
\end{remark}

\section{Some applications}\label{sec.semimartingale}
 We now show that our approach yields simple proofs for some properties of pathwise quadratic variation, which turn out to be useful in the study of pathwise approaches to Ito calculus.
 
Denote $\D:=\D([0,\infty),\R)$ and $\D^{d\times d}:=\D([0,\infty),\R^{d\times d})$ to be the Skorokhod spaces, each of which equipped with a \emph{complete} metric $\delta$ which induces the corresponding Skorokhod (a.k.a J$_1$) topology. Denote $\F$ to be the J$_1$ Borel sigma algebra of $\D$ (a.k.a the canonical sigma algebra generated by coordinates). 
    Recall that   $\Q_{0}^{\pi}$ is the set of paths with finite quadratic variation along $\pi$ in the sense of Def.~\ref{def:2.1}.\\
	
	We now give a criterion for $x\in\D$ to have finite quadratic variation without any reference to the Lebesgue decomposition (\ref{eq:2.2}) on the limit measure $\mu$:
\begin{property}
$x\in Q_{0}^{\pi}$ if and only if $(q_n)$ defined by \eqref{eq:2.6} is a Cauchy sequence in $(\D,\delta)$.
\end{property}

\begin{proof}
This is a consequence of Thm.~\ref{thm:2.7} and that $(\D,\delta)$ is complete.
\end{proof}

	One of the main advantages of having convergence in the J$_1$ topology is that it ensures convergence of jumps in a regulated manner. It comes in handy when accessing the limit of $q_n(t_n)$ as $n\rightarrow\infty$.\\
    
    \begin{property}
    Let $x\in Q_{0}^{\pi}$, for each $t\geq 0$, we define $t'_n:=\max\{t_i<t|t_i\in\pi_n\}$, then

\begin{alignat*}{2}
&t_n\longrightarrow t; t_n\leq t'_n\Longrightarrow q_n({t_n-})&\longrightarrow& [x]({t-}),\\
&t_n\longrightarrow t; t_n<t'_n\Longrightarrow q_n({t_n})&\longrightarrow& [x]({t-}),\\
&t_n\longrightarrow t; t_n\geq t'_n\Longrightarrow q_n({t_n})&\longrightarrow& [x](t),\\
&t_n\longrightarrow t; t_n>t'_n\Longrightarrow q_n({t_n-})&\longrightarrow& [x]({t}).\\
\end{alignat*}
	
    In particular, the sequence $(t'_n)$ is asymptotically unique in the sense that any other sequence $(t''_n)$ meeting the above properties coincides with $(t'_n)$ for  $n$ sufficiently large.
\end{property}

\begin{proof}
This is a consequence of Thm.~\ref{thm:2.7} and \cite[VI.2.1]{jacodshiryaev}.
\end{proof}

	Given a \cadlag process $X$ (i.e. a $(\D,\F)$-measurable random variable), a natural quantity to consider is $\P(X\in\Q_{0}^{\pi})$. This only makes sense however if $\Q_{0}^{\pi}$ is $\F$-measurable. This 'natural' property, not easy to show using the original definition (Def. \ref{def:2.1}), becomes simple thanks to Theorem.~\ref{thm:2.7}:    
\begin{property}[Measurability of $\Q_{0}^{\pi}$]
$\Q_{0}^{\pi}$ is $\F$-measurable.
\end{property}

\begin{proof}
By Thm.~\ref{thm:2.7}, $\Q_{0}^{\pi}=\Q^{\pi}$ and by definition, $\Q^{\pi}$ is the J$_1$ convergence set of\begin{eqnarray*}x\longmapsto\sum_{\pi_{n}\ni t_i\leq\cdot}(x(t_{i+1})-x(t_i))^2,\end{eqnarray*}$n\geq 1$ on $\D$. Since $\D$ is completely metrisable, the claim follows  from \cite[V.3]{doob}.
\end{proof}

    F\"{o}llmer introduced in \cite{follmer1981,follmer1981b} the class of processes with finite quadratic variation ('processus \`a variation quadratique'), defined as  \cadlag processes  such that the sequence \begin{eqnarray*}S_n(t):=\sum_{\pi_{n}\ni t_i\leq t}(X(t_{i+1})-X(t_i))^2 
\end{eqnarray*} converges in probability for every $t$ to an increasing process $[X]$ with Lebesgue decomposition: \begin{eqnarray*}
[X](t)=[X]^{c}(t)+\sum_{s\leq t}\Delta X(s)^t\Delta X(s).
\end{eqnarray*} The pathwise \ito formula \eqref{eq.ito} can be applied to this class of processes, which is strictly larger than the class of semimartingales \cite{coquet2003}.

\begin{theorem}\label{thm:4.4}
Let $\mathbf{X}$ be an $\mathbb{R}^{d}$-valued \cadlag process, define a sequence of $(\D^{d\times d},\delta)$-valued random variables $(\mathbf{q}_n)$ by\begin{eqnarray*}\label{sum}
\mathbf{q}_n(t):=\sum_{\pi_{n}\ni t_i\leq t}(\mathbf{X}(t_{i+1})-\mathbf{X}(t_i))(\mathbf{X}(t_{i+1})-\mathbf{X}(t_i))^{T}
\end{eqnarray*}

then the following properties are equivalent:\begin{enumerate}
\item[(i)] $\mathbf{X}$ is a  process with finite quadratic variation.
\item[(ii)] $(\mathbf{q}_n)$ converges in probability.
\item[(iii)] $(\mathbf{q}_n)$ is a Cauchy sequence in probability.
\end{enumerate}
 In addition, \begin{enumerate}
\item[iv] If $(\mathbf{q}_n)$ converges in probability, the limit is $[\mathbf{X}]$. 
\item[v] The convergence of $(\mathbf{q}_n)$ to $[X]$ is UCP if and only if $\mathbf{X}$ is a continuous process of quadratic variation $[\mathbf{X}]$.
\item[vi] $(\mathbf{q}_n)$ converges (resp. is a Cauchy sequence) in probability if and only if each component sequence of $(\mathbf{q}_n)$ converges (resp. is a Cauchy sequence) in probability.
\end{enumerate}  
\end{theorem}

\begin{proof}
We first remark that $(\D^{d\times d},\delta)$ is a complete separable metric space \cite{jacodshiryaev}, hence by \cite[Lemma  9.2.4]{dudley}, the Cauchy property is equivalent to convergence in probability. By \cite[Thm. 9.2.1]{dudley}, we can pass to subsequences and apply Prop.~\ref{prop:3.5}, Thm.~\ref{thm:3.6} \& Cor.~\ref{cor:3.7} pathwise to $\mathbf{X}$, the claims follow.\end{proof}

{\bf Acknowledgements}\ {We thank Rafal Lochowski, Pietro Siorpaes and the referee for useful comments.}

\providecommand{\bysame}{\leavevmode\hbox to3em{\hrulefill}\thinspace}
\providecommand{\MR}{\relax\ifhmode\unskip\space\fi MR }
\providecommand{\MRhref}[2]{%
  \href{http://www.ams.org/mathscinet-getitem?mr=#1}{#2}
}
\providecommand{\href}[2]{#2}

\end{document}